\documentclass[11pt]{article}
\usepackage{epsfig}
\usepackage{amssymb}
\usepackage{amscd}
\usepackage[all]{xy}
\usepackage{epsf}
\usepackage{amsfonts,amssymb,amsthm,amsmath}

\input cyracc.def

\long\def\comment#1\endcomment{\relax}

\newcounter{subsubsubsection}
\newcounter{subsubsubsubsection}

\makeatletter \@addtoreset{subsubsubsection}{subsubsection}
\@addtoreset{subsubsubsubsection}{subsubsubsection}

\makeatother

\DeclareMathSizes{11.1}{10}{8}{6}

\DeclareMathOperator{\Hom}{Hom}

\newtheorem*{theorem*}{Theorem}

\newtheorem*{conjecture}{Conjecture}

\newtheorem*{lemma}{Lemma}

\newtheorem*{proposition}{Proposition}

\theoremstyle{remark}

\newtheorem*{remark}{Remark}

\newtheorem*{example}{Example}

\theoremstyle{definition}

\newcommand{\eqto}{\mathrel{\stackrel{\sim}{\to}}}

\DeclareMathOperator{\Alt}{Alt}

\newcommand{\hdot}{{\:\protect\raisebox{3pt}{\text{\protect\circle*{1.5}}}}}

\newcommand{\mb}{\hdot}

\newcommand{\Ext}{\mathrm{Ext}}

\newcommand{\g}{\mathfrak{g}}
\newcommand{\Sym}{\mathrm{Sym}}

\newcommand{\K}{\mathrm{K}}

\newcommand{\Der}{\mathrm{Der}}
\newcommand{\Inn}{\mathrm{Inn}}

\newcommand{\xxx}{\otimes\mathbb{C}[[\hbar]]}

\newcommand{\Hoch}{\mathrm{Hoch}}
\newcommand{\poly}{{poly}}

\title{ {\huge  Kontsevich formality and PBW algebras}}
\author{{\LARGE Boris Shoikhet}}
\date{}

\begin{document}\maketitle

{\font\tcyr=wncyi10

\tcyr\cyracc

\font\scyr=wncyr10  \scyr\cyracc \comment \hbox
to\textwidth{\hfil\parbox{80mm}{\tcyr{Zachem ya algebry ne znal podi
po{\u\i}mi se{\u\i}chas}\\{\tcyr Tako{\u\i} polezny{\u\i}
delovo{\u\i} predmet...}}} \vspace{2mm}
\endcomment

\comment \scyr\cyracc  \hbox
to\textwidth{\hfil\parbox{50mm}{\tcyr{Ya vchera fantaziroval,}\\
\tcyr{i fantaziya  mne udalas\cprime...}}} \vspace{2mm}

\hbox to\textwidth{\hfil\parbox{30mm}{\scyr{M.Shcherbakov}}}
\endcomment

\scyr\cyracc  \hbox
to\textwidth{\hfil\parbox{85mm}{\tcyr{Ty po\u\i, moya radost\cprime, ty po\u\i, mo\"e chudo, gitara,}\\
\tcyr{Prozrachnye struny slucha{\u\i}no schastlivogo mira...}}}
\vspace{2mm}

\hbox to\textwidth{\hfil\parbox{18mm}{\scyr{I. Bely\u\i}}}}

\vspace{1cm}

\begin{abstract}
Let $\alpha$ be a polynomial Poisson bivector on a
finite-dimensional vector space $V$ over $\mathbb{C}$. Then
Kontsevich [K97] gives a formula for a quantization $f\star g$ of
the algebra $S(V)^*$. We give a construction of an algebra with the
PBW property defined from $\alpha$ by generators and relations.
Namely, we define an algebra as the quotient of the free tensor
algebra $T(V^*)$ by relations $x_i\otimes x_j-x_j\otimes
x_i=R_{ij}(\hbar)$ where $R_{ij}(\hbar)\in T(V^*)\otimes\hbar
\mathbb{C}[[\hbar]]$, $R_{ij}=\hbar
\Sym(\alpha_{ij})+\mathcal{O}(\hbar^2)$, with one relation for each
pair of $i,j=1...\dim V$. We prove that the constructed algebra
obeys the PBW property, and this is a generalization of the
Poincar\'{e}-Birkhoff-Witt theorem. In the case of a linear Poisson
structure we get a new very conceptual proof of the PBW theorem
itself, and for a quadratic Poisson structure we get an object
closely related to a quantum $R$-matrix on $V$.

The construction uses the Kontsevich formality. Namely, our
quantities $R_{ij}(\hbar)\in T(V^*)\otimes \mathbb{C}[[\hbar]]$ are
written directly in Kontsevich integrals from [K97], but in a sense
of dual graphs than the graphs used in the deformation quantization.
We conjecture that the relation $x_i\otimes x_j-x_j\otimes
x_i=R_{ij}(\hbar)$ holds in the Kontsevich star-algebra, when we
replace $\otimes$ by $\star$. This conjecture implies in particular
that our algebra is isomorphic to the Kontsevich star-algebra with
the same $\alpha$, but also gives a highly-nontrivial identity on
Kontsevich integrals. Probably, it is a particular case of a more
general duality acting on the AKSZ model on open disc, used by
Kontsevich in his proof of the formality conjecture.
\end{abstract}

\section*{Introduction}
\subsection{}
Let $\alpha$ be a polynomial Poisson bivector on a vector space $V$.
There are two ways how one can think about what is to quantize
$\alpha$.

The first way, called deformation quantization, is well-known. One
looks for an associative product on the space $C^\infty(V)$ of
smooth functions on $V$ having a form
\begin{equation}\label{intro1}
f\star g=f\cdot g+\hbar \frac12\{f,g\}_\alpha+\hbar^2
B_2(f,g)+\hbar^3 B_3(f,g)+\dots
\end{equation}
where $\{f,g\}_\alpha=\alpha(df\wedge dg)$ is the corresponding
Poisson bracket, and $B_i\colon C^\infty(V)^{\otimes 2}\to
C^\infty(V)$ are some local maps, which in practice are
bi-differential operators. The associativity condition
\begin{equation}\label{intro2}
(f\star g)\star h=f\star(g\star h)
\end{equation}
for any $f,g,h\in C^\infty(V)$ is then an infinite sequence of
quadratic identities on $B_i$'s. This problem was solved by
M.Kontsevich in 1997, and the solution uses a two-dimensional
topological quantum field theory on open disc.

There is another way how to think about the "quantization". Namely,
we are looking for an algebra which is a quotient of the tensor
algebra $T(V^*)$ by relations of the form
\begin{equation}\label{intro3}
x_i\otimes x_j-x_j\otimes x_i=\hbar\Sym(\alpha_{ij})+\hbar^2
\omega_2+\hbar^3\omega_3+\dots
\end{equation}
Here $\{x_i\}$ is a basis in $V^*$, and
$\alpha=\sum_{i,j}\alpha_{ij}\partial_i\wedge\partial_j$, where
$\alpha_{ij}$ are polynomials. The symmetrization
$\Sym(\alpha_{ij})$ is the element of $T(x_1,\dots, x_n)$ which is
given by the full symmetrization, $\Sym\colon
\mathbb{C}[x_1,\dots,x_n]\to T(x_1,\dots, x_n)$. Here
$\omega_2,\omega_3,\dots$ are some elements of $T(x_1,\dots, x_n)$,
they depend on the pair $(i,j)$ of indices. So, the problem is to
find such $\omega_k(i,j)$'s such that the obtained algebra obeys the
Poincar\'{e}-Birkhoff-Witt (PBW) property. Let us formulate this
property.

Denote by $A_\hbar$ the algebra over $\mathbb{C}[[\hbar]]$ given by
(\ref{intro3}). We have the following decreasing filtration of
$A_\hbar$:

\begin{equation}\label{intro4}
A_\hbar\supset \hbar A_\hbar \supset\hbar^2A_\hbar\supset\dots
\end{equation}
It is an algebra filtration, that is, $(\hbar^i A_\hbar)\cdot
(\hbar^j A_\hbar)\subset\hbar^{i+j}A_\hbar$. Consider the associated
graded algebra $\mathrm{gr}A_\hbar=\oplus_{k\ge 0}\hbar^k
A_\hbar/\hbar^{k+1}A_\hbar$. The PBW condition then is that there is
an $\hbar$-linear isomorphism $\mathrm{gr}A_\hbar\simeq
S(V^*)\otimes \mathbb{C}[[\hbar]]$ with which $\hbar^k
A_\hbar/\hbar^{k+1}A_\hbar\simeq\hbar^k S(V^*)$.

If the algebra $A_\hbar$ is given by (\ref{intro3}), the condition
above always holds for $k=0$, and the only condition which is
necessary for $k=1$ is $\alpha_{ij}=-\alpha_{ji}$. For $k\ge 2$ in
general the component $\hbar^k A_\hbar/\hbar^{k+1}A_\hbar$ is {\it
less} than $\hbar^k S(V^*)$. We can go one step further and find the
necessary condition for $k=2$.

Consider $\Alt_{i,j,k}[x_i,[x_j,x_k]]$ where $[a,b]=a\star b -b\star
a$. This expression is 0 for any associative algebra. Actually this
condition is an infinite sequence of conditions, but its image in
$\hbar^2 A_\hbar/\hbar^3 A_\hbar$ depends only on $\{\alpha_{ij}\}$.
It is exactly the condition $\{\alpha, \alpha\}=0$ where the bracket
is the Schouten-Nijenhuis bracket. If this condition is not
satisfied, we get a  nonzero element in $\hbar^2 S(V^*)$ which is
zero in $\hbar^2A_\hbar/\hbar^3A_\hbar$. One can check that if the
Schouten-Nijenhuis bracket $\{\alpha,\alpha\}$ is zero, then
$\hbar^2 A_\hbar/\hbar^3A_\hbar\simeq \hbar^2 S(V^*)$.

For higher orders in $\hbar$ we have more complicated conditions.
The claim is that if $\alpha$ is a Poisson bivector, we can find all
$\omega_k(i,j)$ in (\ref{intro3}) such that the algebra defined by
(\ref{intro3}) is a PBW algebra. The reader can see from this
discussion that this question is very close to the classical
question of deformation quantization.

\subsection{}
The Kontsevich deformation quantization formula gives us a PBW
algebra associated with a Poisson bivector $\alpha$. Indeed, we have
some formula like $x_i\star x_j=x_i\cdot x_j
+\frac12\hbar\alpha_{ij}(x_1,\dots, x_n)+\hbar^2(\dots)+\dots$. Here
$\star$ is the Kontsevich star product. We see that $x_i\star
x_j-x_j\star x_i$ starts with the first order in $\hbar$. The
right-hand side is an element in
$S(V^*)\otimes\hbar\mathbb{C}[[\hbar]]$. If we express iteratively
the right-hand side as sum of monomials of the form $x_{i_1}\star
x_{i_2}\star\dots\star x_{i_k}$, we can then replace $\star$ by
$\otimes$ and get a PBW algebra. It will be indeed a PBW algebra,
because the associated graded algebra has at most the size as
$S(V^*)\otimes\mathbb{C}[[\hbar]]$, and it can not have a less size
because this relation holds in the Kontsevich algebra. This proves
in particular that one doesn't need any other relation to define the
Kontsevich star-algebra.

A lack of this construction is that we apply Kontsevich formula, or
rather something like a reverse to it, infinitely many times. The
coefficients will depend on the Kontsevich integrals in deformation
quantization, but actually will be much more complicated.

Our solution uses the Kontsevich integrals of, in a sense, dual
admissible graphs. In particular, they are given directly in
Kontsevich integrals, without any iterative process, but of dual
graphs. We conjecture that our relation holds exactly in the
Kontsevich star-algebra. If this conjecture is true, we get a very
complicated relation between Kontsevich integrals.

\subsection{}
It is clear that if the conjecture described in the previous
Subsection is true, it should have an analogue for all integrals,
not only for particular graphs involving in the deformation
quantization formula. We say that we want to lift this conjecture
"on the level of complexes". Moreover, we believe that some "Koszul
duality" acts on the entire AKSZ model on open disc, and we want to
express this duality mathematically. We are going to consider this
question in the sequel.
\subsection{}
All quadratic PBW algebras are Koszul. One easily sees that if one
starts with a quadratic Poisson bivector $\alpha$, the relations we
get are quadratic, namely, all $\omega_k(i,j)$ are elements in
$V^*\otimes V^*\subset T(V^*)$. The dg Lie algebra of polyvector
fiels on $V$ is isomorphic to the dg Lie algebra of polyvector
fields on $V^*[1]$, and this isomorphism preserves quadratic
bivector fields. Then, our constructions give two PBW algebras from
a quadratic Poisson bivector on $V$: one is the quotient of $T(V^*)$
by some quadratic relations, and another is a quotient of $T(V[-1])$
by some other quadratic relations. It is natural to conjecture that
the two algebras are Koszul dual. This explains our notation "Koszul
duality in deformation quantization". We are going to consider these
questions in a sequel paper.

\subsection{}
We would like especially to note, that the $L_\infty$ map
\begin{equation}\label{intro100}
\Theta\colon T_\poly(V)\to
\Der(CoBar^\mb(\Lambda^-(V^*)))/\Inn(CoBar^\mb(\Lambda^-(V^*)))
\end{equation}
which is the composition of the map $\Phi$ from Section 1.4 with the
Kontsevich formality $L_\infty$ map, can be considered as a
"non-commutative analogue" of the usual map (actually an
isomorphism) $W_n\to \Der(S(V^*))$ where $W_n$ is the Lie algebra of
polynomial vector fields on $V$. Probably it is possible to define
some "quasi-manifolds" by replacing the local coordinate ring
$S(V^*)$ by its free resolution $CoBar^\mb(\Lambda^-(V^*))$, and by
replacing the usual transition functions by $A_\infty$
quasi-isomorphisms with the natural compatibility property. We need
to invert some of these $A_\infty$ quasi-isomorphisms in the
compatibility equation on the "triple intersection" $U_i\cap U_j\cap
U_k$, therefore, we should work in the Quillen homotopical category.
Any usual manifold gives us a trivial example of a quasi-manifold by
means of the canonical projection of the resolution to the algebra.
The corresponding Lie algebra of infinitesimal symmetries in the
sense of the formal geometry will be then $T_\poly(V)$, by means of
the quasi-isomorphism (\ref{intro100}).

A first difficulty in the realization of this program is that our
map (\ref{intro100}), and even the resolution, are defined for
polynomial algebras, not algebras of smooth functions. Nevertheless,
the author is sure that an algebraic analog of this construction,
leading to a rigorous definition of a "quasi-manifold", should
exist. We hope to clarify it in the sequel.

\subsection*{Acknowledgements}
I am thankful to Misha Bershtein, Pavel Etingof, Giovanni Felder,
and especially to Borya Feigin, Victor Ginzburg and to Bernhard
Keller for interesting and useful discussions. I am grateful to Jim
Stasheff for many remarks and suggestions. The work was partially
supported by the research grant R1F105L15 of the University of
Luxembourg.

\section{The main construction}
\subsection{}
First of all, recall here the Stasheff's definition of the
Hochschild cohomological complex of an associative algebra $A$.

Consider the shifted vector space $W=A[1]$, and the cofree
coassociative coalgebra $C(W)$ (co)generated by $W$. As a graded
vector space, $C(W)=T(A[1])$, the free tensor space. The coproduct
is:

\begin{equation}\label{eq1}
\begin{aligned}
\ &\Delta(a_1\otimes a_2\otimes\dots\otimes a_k)=\sum_{i=1}^{k-1}
(a_1\otimes\dots \otimes a_i)\bigotimes (a_{i+1}\otimes\dots\otimes
a_k)
\end{aligned}
\end{equation}

Consider the Lie algebra $CoDer(C(A[1]))$ of all coderivations of
this coalgebra. As the coalgebra is free, any coderivation $D$ (if
it is graded) is uniquely defined by a map $\Psi_D\colon A^{\otimes
k}\to A$, and the degree of this coderivation is $k-1$ (in
conditions that $A$ is not graded). The bracket
$[\Psi_{D_1},\Psi_{D_2}]$ is again a coderivation. Define the
Hochschild Lie algebra as $\Hoch^\mb(A)=CoDer^\mb(C(A[1]))$. To
define the complex structure on it, consider the particular
coderivation $D_m$ of degree +1 from the product $m\colon A^{\otimes
2}\to A$, which is the product in the associative algebra $A$. The
condition $[D_m,D_m]=0$ is equivalent to the associativity of $m$.
Define the differential on $CoDer^\mb(C(A[1]))$ as
$d(\Psi)=[D_m,\Psi]$. In this way we get a dg Lie algebra. The
differential is called the Hochschild differential, and the bracket
is called the Gerstenhaber bracket. The definition of these
structures given here is due to J.Stasheff.

\subsection{The explicit definition}
Here we relate the Stasheff's definition of the Hochschild
cohomological complex with the usual one.

The concept of a coderivation of a (co)free coalgebra is dual to the
concept of a derivation of a free algebra. Let $L$ be a vector
space, and let $T(L)$ be the free tensor algebra generated by the
vector space $L$. Let $D\colon T(L)\to T(L)$ be a derivation, then
it is uniquely defined by its value $D_L\colon L\to T(L)$ on the
generators, and any $D_L$ defines a derivation $D$ of the free
algebra $T(L)$. If we would like to consider only graded
derivations, we restrict ourselves by the maps $D_L\colon L\to
L^{\otimes k}$ for $k\ge 0$.

Dually, a coderivation $D$ of the cofree coalgebra $C(P)$
cogenerated by a vector space $P$ is uniquely defined by the
restriction to cogenerators, that is, by a map $D_P\colon C(P)\to
P$, or, if we consider the graded coderivations, the map $D_P$ is a
map $D_P\colon P^{\otimes k}\to P$ for $k\ge 0$.

In our case of the definition of the cohomological Hochschild
complex of an associative algebra $A$, we have $P=A[1]$. Then the
coderivations of the grading $k$ form the vector space
$\Hoch^k(A)=\Hom(A^{\otimes (k+1)},A)$, $k\ge -1$. Now we can deduce
the differential and the Gerstenhaber bracket from the Stasheff's
construction. The answer is the following:

For $\Psi\in\Hom(A^{\otimes k},A)$ the cochain $d\Psi\in
\Hom(A^{\otimes (k+1)},A)$ is given by the formula:

\begin{equation}\label{eq20}
\begin{aligned}
\ &d\Psi(a_0\otimes\dots\otimes a_k)=a_0\Psi(a_1\otimes\dots\otimes a_k)+\\
&+\sum_{i=0}^{k-1}(-1)^{i+1}\Psi(a_0\otimes\dots\otimes
a_{i-1}\otimes(a_ia_{i+1})\otimes a_{i+2}\otimes\dots\otimes a_k)+\\
&+(-1)^{k+1}\Psi(a_0\otimes\dots\otimes a_{k-1})a_k
\end{aligned}
\end{equation}

For $\Psi_1\in\Hom(A^{\otimes(k+1)},A)$ and
$\Psi_2\in\Hom(A^{\otimes (l+1)},A)$ the bracket
$[\Psi_1,\Psi_2]=\Psi_1\circ\Psi_2-(-1)^{kl}\Psi_2\circ\Psi_1$ where

\begin{equation}\label{eq21}
\begin{aligned}
\ & (\Psi_1\circ\Psi_2)(a_0\otimes\dots\otimes a_{k+l})=\\
&\sum_{i=0}^k(-1)^{il}\Psi_1(a_0\otimes\dots\otimes
a_{i-1}\otimes\Psi_2(a_i\otimes\dots\otimes a_{i+l})\otimes
a_{i+l+1}\otimes\dots\otimes a_{k+l})
\end{aligned}
\end{equation}

\subsection{The (co)bar-complex}
Here we recall the definition of the (co)bar-complex of an
associative (co)algebra. When the (co)algebra contains (co)unit, the
(co)bar-complex is acyclic, and when the (co)algebra is the kernel
of the augmentation of a quadratic Koszul algebra, this concept is
closely related to the Koszul duality.

Let $A$ be an associative algebra. Then its bar-complex is
$$
\dots\rightarrow A^{\otimes 3}\rightarrow A^{\otimes 2}\rightarrow
A\rightarrow 0
$$
where $\deg A^{\otimes k}=-k+1$, and the differential $d\colon
A^{\otimes k}\to A^{\otimes (k-1)}$ is given as follows:

\begin{equation}\label{eq23}
d(a_1\otimes\dots\otimes a_k)=(a_1a_2)\otimes a_3\otimes\dots\otimes
a_k-a_1\otimes (a_2a_3)\otimes\dots\otimes a_k+\dots
+(-1)^{k}a_1\otimes\dots\otimes a_{k-2}\otimes (a_{k-1}a_k)
\end{equation}
If the algebra $A$ has unit, the bar-complex of $A$ is acyclic in
all degrees. Indeed, the map
$$
a_1\otimes \dots\otimes a_k\mapsto 1\otimes a_1\otimes \dots\otimes
a_k
$$
is a contracting homotopy.

Suppose now that the algebra $A$ does not contain unit, and $A=B^+$
is the kernel of an augmentation map $\varepsilon\colon
B\to\mathbb{C}$. (The map $\varepsilon$ is a surjective map of
algebras, in particular, it maps $1$ to $1$). Then the cohomology of
the bar-complex of $A$ is equal to the dual space
$\Ext^\mb_{B-Mod}(\mathbb{C},\mathbb{C})$.

Indeed, for any $B$-module $M$, we have the following free
resolution of $M$:

\begin{equation}\label{eq25}
\dots B\otimes \overline{B}\otimes\overline{B}\otimes M\rightarrow
B\otimes\overline{B}\otimes M\rightarrow B\otimes M\rightarrow
M\rightarrow 0
\end{equation}
with the differential analogous to the bar-differential.

Consider the case $M=\mathbb{C}$. We can compute
$\Ext_{B-Mod}^\mb(\mathbb{C},\mathbb{C})$ using this resolution. In
the answer we get the cohomology of the complex dual to the
bar-complex of $A=B^+$.

The complex dual to the bar-complex of $A$ is the {\it
cobar-complex} for the coalgebra $A^*$. This cobar-complex is an
associative dg algebra, and it is a free algebra, which by previous
is a free resolution of the algebra
$\Ext_{B-Mod}(\mathbb{C},\mathbb{C})$. For the sequel we write down
explicitly the cobar-complex of a coassociative coalgebra $Q$:

\begin{equation}\label{eq26}
0\rightarrow Q\rightarrow Q\otimes Q\rightarrow Q\otimes Q\otimes
Q\rightarrow\dots
\end{equation}
and the differential $\delta Q^{\otimes k}\to Q^{\otimes (k+1)}$ is
\begin{equation}\label{eq27}
\delta(q_1\otimes\dots\otimes q_k)=(\Delta q_1)\otimes
q_2\otimes\dots\otimes q_k-q_1\otimes(\Delta q_2)\otimes
\dots\otimes q_k+\dots+(-1)^{k-1}q_1\otimes \dots\otimes
q_{k-1}\otimes (\Delta q_k)
\end{equation}
where $\Delta\colon Q\to Q^{\otimes 2}$ is the coproduct.

In the case when $B=S(V)$ is the symmetric algebra,
$\Ext_{B-Mod}(\mathbb{C},\mathbb{C})$ is the exterior algebra
$\Lambda (V^*)=S(V^*[-1])$, and vise versa. In this way, we get a
free resolution of the symmetric (exterior) algebra.

\begin{example}
Here we construct explicitly the free cobar-resolution
$\mathcal{R}^\mb$ of the algebra $\mathbb{C}[x_1,x_2]$ of
polynomials on two variables. As a graded algebra, $\mathcal{R}^\mb$
is the free algebra
$\mathcal{R}^\mb=\mathrm{Free}(x_1,x_2,\xi_{12})$ where $\deg
x_1=\deg x_2=0$, $\deg\xi_{12}=-1$. The differential is $0$ on
$x_1,x_2$, $d(\xi_{12})=x_1\otimes x_2-x_2\otimes x_1$, and
satisfies the graded Leibniz rule. In degree 0 we have the tensor
algebra $T(x_1,x_2)$, differential is 0 on degree 0 (there are no
elements in degree 1). In degree -1, a general element is a
non-commutative word in $x_1,x_2,\xi_{12}$ in which $\xi_{12}$
occurs exactly one time. For example, it could be a word $x_2\otimes
x_1\otimes x_2\otimes \xi_{12}\otimes x_1\otimes x_2$. The image of
the differential is then exactly the two-sided ideal in the tensor
algebra $T(x_1,x_2)$ generated by $x_1\otimes x_2-x_2\otimes x_1$.
Then, the 0-th cohomology is $\mathbb{C}[x_1,x_2]$. It follows from
the discussion above that all higher cohomology is 0.
\end{example}

\begin{example}
Consider the case of the algebra $\mathbb{C}[x_1,x_2,x_3]$. Then
$\mathcal{R}^\mb=\mathrm{Free}(x_1,x_2,x_3,\xi_{12},
\xi_{23},\xi_{13},\xi_{123})$ with $\deg x_i=0$, $\deg \xi_{ij}=-1$,
$\deg\xi_{123}=-2$. The differential is 0 on $x_1,x_2,x_3$,
$d(\xi_{ij})=x_i\otimes x_j-x_j\otimes x_i$, and
$d(\xi_{123})=(x_1\otimes \xi_{23}+x_2\otimes \xi_{31}+x_3\otimes
\xi_{12})+(\xi_{23}\otimes x_1+\xi_{31}\otimes x_2+\xi_{12}\otimes
x_3)$. Here we set $\xi_{ij}=-\xi_{ji}$. Then cohomology in degree 0
is equal to the quotient of the free algebra $T(x_1,x_2,x_3)$ by the
two-sided ideal generated by $x_i\otimes x_j-x_j\otimes x_i$, that
is, the algebra $\mathbb{C}[x_1,x_2,x_3]$. It follows from our
general discussion in this Subsection that the higher cohomology
vanishes.
\end{example}

\subsection{The main construction} Here we construct a quasi-isomorphic map of dg Lie
algebras $\Phi\colon \Hoch^\mb(S(V))/\mathbb{C}\to
\Der(CoBar^\mb(S(V^*)^+))/\Inn (CoBar^\mb(S(V^*)^+))$.

Let $\Psi\in\Hom ((S(V))^{\otimes k},S(V))$ be a $k$-cochain. Denote
$V^*=W$, then we can consider $\Psi$ the corresponding cochain in
$\Hom(S(W),(S(W))^{\otimes k})$. Here we consider $S(W)$ as {\it
coalgebra}. Then this cochain may be considered as a derivation in
$\Der(CoBar^\mb(S(W)))$. We would like to attach to it a derivation
in $\Der(CoBar^\mb(S(W)^+))$, maybe modulo an inner derivation. So,
we would like to show that there exist a map $\Phi\colon
\Der(CoBar^\mb(S(W)))\to\Der(CoBar^\mb(S(W)^+))$ such that the
diagram

\begin{equation}\label{diagram1}
\xymatrix{ \Der(CoBar^\mb(S(W)))\ar[r]^{\delta}  \ar[d]^{\Phi}&
\Der(CoBar^\mb(S(W)))\ar[d]^{\Phi}\\
\Der(CoBar^\mb(S(W)^+))\ar[r]^{\delta} &\Der(CoBar^\mb(S(W)^+))}
\end{equation}
is commutative {\it modulo inner derivations} (here $\delta$ is the
cobar differential).

In the coalgebra $S(W)$ the coproduct is given by the formula
\begin{equation}\label{eq27}
\Delta(x_1\dots x_k)=1\otimes (x_1\dots x_k)+\sum_{\{i_1\dots
i_a\}\sqcup\{j_1\dots j_b\}=\{1\dots k\}}(x_{i_1}\dots
x_{i_a})\otimes(x_{j_1}\dots x_{j_b})+(x_1\dots x_k)\otimes 1
\end{equation}
and in the coalgebra $S(W)^+$ the coproduct is given by the same
formula without the first and the last summands, which contain 1's.

Therefore the projection $p\colon S(W)\to S(W)^+$ is a map of
coalgebras (dual to the imbedding of algebras), and the imbedding
$i\colon S(W)^+\to S(W)$ is {\it not}.

If $\Psi\colon S(W)\to S(W)^{\otimes k}$ is as above, we define
$(\Phi(\Psi))(\sigma)=p^{\otimes k}(\Psi(i(\sigma)))\in\Hom (S(W)^+,
(S(W)^+)^{\otimes k})$.

Now we check the commutativity of the diagram (\ref{diagram1})
modulo inner derivations. It is clear that

\begin{equation}\label{eq28}
(\Phi\circ\delta)(\sigma)-(\delta\circ\Phi)(\sigma)=p^{\otimes
k}(\Psi(1))\otimes\sigma\pm \sigma\otimes p^{\otimes k}(\Psi(1))
\end{equation}
which is an inner derivation $ad(p^{\otimes k}(\Psi(1)))$.

We have defined a map $\Phi_1\colon \Der(CoBar^\mb(S(W)))\to
\Der(CoBar^\mb(S(W)^+))/\Inn(CoBar^\mb(S(W)^+))$. The first dg Lie
algebra is clearly isomorphic to the Hochschild cohomological
complex of the algebra $S(V)$ modulo constants, and we can consider
the map $\Phi_1$ as a map
$$
\Phi_1\colon\Hoch^\mb(S(V))/\mathbb{C}\to\Der(CoBar^\mb(S(V^*)^+))/\Inn(CoBar^\mb(S(V^*)^+))
$$
Let us note that the cobar complex $CoBar^\mb(S(V^*)^+)$ is a free
resolution of the Koszul dual algebra $\Lambda(V^*)$.

Now we have the following result:

\begin{proposition}
The map $\Phi_1$ is a quasi-isomorphism of dg Lie algebras.
\end{proposition}
It is clear that $\Phi_1$ is a map of dg Lie algebras, one only
needs to proof that it is a quasi-isomorphism of complexes. We will
not use this result, the proof will appear somewhere. Let us,
however, outline some naive ideas behind the proof in the next
Subsection.

\subsection{} We prove firstly that the cohomology of
$\Der(CoBar^\mb(S(V^*)^+))/\Inn(CoBar^\mb(S(V^*)^+))$ is isomorphic
to the Hochschild cohomology $\Hoch^\mb(CoBar^\mb(S(V^*)^+))$ of the
dg algebra $CoBar^\mb(S(V^*)^+)$. Let us note that it is not true
that the Hochschild cohomology $\Hoch^\mb(CoBar^\mb(B))/\mathbb{C}$
is isomorphic to the cohomology of
$\Der(CoBar^\mb(B))/\Inn(CoBar^\mb(B))$ for any coassociative
coalgebra $B$. To have this property, $B$ should be {\it
cocomplete}, that is

\begin{equation}\label{14proof1}
B=\bigcup_{n\ge 1}\mathrm{Ker}(\Delta^n\colon B\to B^{\otimes n+1})
\end{equation}

This fact is proven in [Lef]. The coalgebras $S^+(V)$,
$\Lambda^-(V)$ are clearly cocomplete for any vector space $V$. The
coalgebra $S(V)$ is not cocomplete, and the property fails for it.
Indeed, let us suppose that
$\Der(CoBar^\mb(S(V)))/\Inn(CoBar^\mb(S(V)))$ has the same
cohomology that $\Hoch^\mb(CoBar^\mb(S(V)))$. The dg algebra
$CoBar^\mb(S(V))$ is acyclic, and the Hochschild cohomology of
quasi-isomorphic algebras are the same; we conclude, that
$HH^\mb(CoBar^\mb(S(V)))=0$. But the cohomology of
$\Der^\mb(CoBar^\mb(S(V)))/\Inn(CoBar^\mb(S(V)))$ is equal to
polyvector fields $T_\poly(V)/\mathbb{C}$. Indeed, consider the
short exact sequence of complexes
\begin{equation}\label{14proof2}
\begin{aligned}
\ &0\rightarrow \Inn(CoBar^\mb(S(V)))\rightarrow
\Der(CoBar^\mb(S(V)))\rightarrow\\
&\rightarrow \Der(CoBar^\mb(S(V)))/\Inn(CoBar^\mb(S(V)))\rightarrow
0
\end{aligned}
\end{equation}
As algebra, $CoBar^\mb(S(V))$ is free, therefore,
$\Inn(CoBar^\mb(S(V)))=CoBar^\mb(S(V))/\mathbb{C}$ is acyclic. From
the long exact sequence associated with (\ref{14proof2}) one has
$H^\mb(\Der(CoBar^\mb(S(V)))/\Inn(CoBar^\mb(S(V))))=H^\mb(CoBar^\mb(S(V)))$,
and the latter is Hochschild cohomology modulo constants of the
algebra $S(V^*)$ by the Stasheff's construction (see Section 1.1).
This cohomology is equal to polyvector fields on $V$ by the
Hochschild-Kostant-Rosenberg theorem.

This example was surprising for the author, because at first look
one has the following "proof" of the statement for general, not only
cocomplete, coalgebra $B$. Consider the Hochschild complex of
$CoBar^\mb(B)$ as a bicomplex. It has two differentials: the
horizontal one is the cobar-differential, and the vertical one is
the Hochschild differential. The bicomplex is placed in the I and II
quarters. The spectral sequence which firstly computes the vertical
(Hochschild) differential converges. Compute firstly the cohomology
with respect to the Hochschild (vertical) differential. As the
coalgebra $CoBar^\mb(B)$ is free, it has Hochschild cohomology only
in degrees 0 and 1; in degree 0 it is $\mathbb{C}$, and in degree 1
it is $\Der(CoBar^\mb(B))/\Inn(CoBar^\mb(B))$. Only the terms
$E_1^{0,0}$ and $E_1^{1,q}$, $q\ge 0$ are nonzero. The spectral
sequence collapses at the term $E_2$. This completes the "proof".

A solution to this contradiction was found in discussions with
Bernhard Keller. If $A^\mb$ is a dg associative algebra with
infinitely many nonzero grading components, there are two possible
definitions of the Hochschild cohomological complex of $A^\mb$.
Namely, the Hochschild complex of $A^\mb$ is a bicomplex, and we can
take the {\it sum} total complex and the {\it product} total
complex. Suppose we take the product total complex, then there is
{\it only one} filtration of this product total complex. Namely, it
is the filtration by subcomplexes $F_k(Tot^\mb)=\{K^{p,q}, q\ge
k\}$. This filtration satisfies the condition $\cup_k
F_k(Tot^\mb)=Tot^\mb$. Another one, $G_k(Tot^\mb)=\{K^{p,q}, p\ge
k\}$ is a filtration only on the sum total complex, that is, for the
product total complex it does {\it not} satisfy $\cup_k
G_k(Tot^\mb)=Tot^\mb$. We actually use the filtration $G_k$ to have
the (vertical) Hochschild differential in the term $E_0$. We see
that this speculation fails, because $G_k$ is not a filtration of
the product total complex. Well, then we can take the sum total
complex. But why we are guaranteed that $\Hoch^\mb(Cobar^\mb(S(V)))$
has 0 cohomology? We use the fact that the Hochschild cohomology of
quasi-isomorphic dg algebras are the same, and $CoBar^\mb(S(V))$ is
quasi-isomorphic to 0 bacause the coalgebra $S(V)$ has a counit (see
Section 1.3). But this fact is proven only for the product total
complex definition of the Hochschild cohomological complex.
Moreover, our discussion shows that this fact is not true for the
sum total complex, and our "proof" fails either.

\section{Applications to deformation quantization}
\subsection{A lemma}
We start with the following lemma, which is a formal version of the
semi-continuity of cohomology of a complex depending on a parameter,
which says that in a "singular value" of the parameter the
cohomology may only raise:
\begin{lemma}
\begin{itemize}
Let $\mathcal{R}^\mb$ be a $\mathbb{Z}_{\le 0}$-graded complex with
differential $d_0$, such that $H^i(\mathcal{R}^\mb)$ vanishes for
all $i\ne 0$.
\item[(i)]Consider
$\mathcal{R}^\mb_{\hbar}=\mathcal{R}^\mb\otimes\mathbb{C}[\hbar]$.
Let $d_\hbar\colon \mathcal{R}^\mb_\hbar\to
\hbar\mathcal{R}^{\mb+1}_\hbar$ be an $\hbar$-linear map of degree
+1 $d_\hbar=d_0+\hbar d_1+\hbar^2 d_2+\dots \hbar^n d_n$ (a finite
sum) such that
$$
d_\hbar^2=0
$$
Denote by $H^\mb_\hbar$ the cohomology of the complex
$(\mathcal{R}^\mb_\hbar,d_\hbar)$. Consider the filtration
$\mathcal{R}^\mb_\hbar\supset\hbar\mathcal{R}^\mb_\hbar\supset\hbar^2\mathcal{R}^\mb_\hbar\supset\dots$,
and the induced filtration on $H^\mb_\hbar$: $F_i
H^\mb_\hbar=\mathrm{Im}(i\colon H^\mb(\hbar^i \mathcal{R}^\mb_\hbar,
d_\hbar)\to H^\mb(\mathcal{R}^\mb_\hbar, d_\hbar))$. Then there are
canonical isomorphisms of $\mathbb{C}[\hbar]$-modules
$$
F_i H^0_\hbar/F_{i+1}H^0_\hbar\eqto \hbar^i H^0(\mathcal{R}^\mb,
d_0)
$$
and $F_i H^k_\hbar/F_{i+1}H^k_\hbar=0$ for $k<0$,
\item[(ii)] The same statement for
$\mathcal{R}^\mb_{[[\hbar]]}=\mathcal{R}^\mb\otimes\mathbb{C}[[\hbar]]$,
and $d_\hbar=d_0+\hbar d_1+\hbar^2 d_2+\dots$, possibly an infinite
sum, $d_\hbar^2=0$.
\end{itemize}
\begin{proof}
Before a rigorous proof, let us make a remark. Consider the
filtration
$$
\mathcal{R}^\mb_\hbar\supset\hbar\mathcal{R}^\mb_\hbar\supset\hbar^2\mathcal{R}_\hbar\supset\dots
$$
of the complex $\mathcal{R}_\hbar^\mb$ with the differential
$d_\hbar$. Let us try to compute the cohomology of
$\mathcal{R}_\hbar^\mb$ by the spectral sequence corresponding to
this filtration. The term
$E_0^{p,q}=\hbar^p\mathcal{R}^{p+q}_\hbar/\hbar^{p+1}\mathcal{R}^{p+q}_\hbar$,
and $d_i$, $i>0$ act by 0 on $E^{p,q}_0$. Therefore, the cohomology
in this term is the cohomology of the differential $d_0$ and is
$E_1^{p,-p}=\hbar^p H^0(\mathcal{R}^\mb,d_0)/\hbar^{p+1}
H^0(\mathcal{R}^\mb,d_0)$ and $E_1^{p,q}=0$ for $q\ne -p$.
Therefore, the spectral sequence collapses at the term $E_1$ by the
dimensional reasons. However, we are not guaranteed that the
spectral sequence converges to the associated graded space with
respect to the induced filtration on cohomology: the spectral
sequence "lives" in the III-rd quarter, and it does not follow from
"dimensional reasons". It is not a rigorous proof, but it somehow
explains the situation.

Now we pass to a rigorous proof. We prove the both statements (i)
and (ii) simultaneously.

Consider the short exact sequences of complexes $S_k$:
\begin{equation}\label{lemma21new1}
0\rightarrow \hbar^{k+1}\mathcal{R}^\mb_\hbar\rightarrow
\hbar^k\mathcal{R}^\mb_\hbar\rightarrow
\hbar^{k}\mathcal{R}^\mb_\hbar/\hbar^{k+1}\mathcal{R}^\mb_\hbar\rightarrow
0
\end{equation}
The complex
$\hbar^{k}\mathcal{R}^\mb_\hbar/\hbar^{k+1}\mathcal{R}^\mb_\hbar$
has the differential $d_0$ because all higher differentials vanish.
Therefore, in the long exact sequence in cohomology corresponding to
$S_k$ we have many zero spaces, namely,
$H^\ell(\hbar^{k}\mathcal{R}^\mb_\hbar/\hbar^{k+1}\mathcal{R}^\mb_\hbar)$
for $\ell\le -1$. Then the long exact sequence proves that the
imbedding
$\hbar^{k+1}\mathcal{R}^\mb_\hbar\hookrightarrow\hbar^k\mathcal{R}^\mb_\hbar$
induces an isomorphism on $\ell$-th cohomology for all $\ell\le -1$.
Consider the end fragment of the long exact sequence:
\begin{equation}\label{lemma21new2}
\begin{aligned}
\ &\dots\rightarrow H^1(\hbar^{k+1}\mathcal{R}^\mb_\hbar)\rightarrow
H^1(\hbar^k\mathcal{R}^\mb_\hbar)\rightarrow 0\rightarrow\\
&\rightarrow H^0(\hbar^{k+1}\mathcal{R}^\mb_\hbar)\rightarrow
H^0(\hbar^{k}\mathcal{R}^\mb_\hbar)\rightarrow
H^0(\hbar^{k}\mathcal{R}^\mb_\hbar/\hbar^{k+1}\mathcal{R}^\mb_\hbar)\rightarrow
0
\end{aligned}
\end{equation}
which proves all assertions of the lemma.

Lemma is proven.
\end{proof}
\end{lemma}

\begin{remark}
Contrary with the case of the lemma, consider the case when the
complex $(\mathcal{R}^\mb, d_0)$ is $\mathbb{Z}_{\ge 0}$-graded, and
again only $H^0(\mathcal{R},d_0)\ne 0$. Then the case (i) of lemma
fails, and only (ii) is true. The proof goes as follows:

By the long exact sequence arguments one needs to prove only the
surjectivity of the map $H^0(\hbar^k \mathcal{R}^\mb_{[[\hbar]]})\to
H^0(\hbar^k \mathcal{R}^\mb_{[[\hbar]]}/\hbar^{k+1}
\mathcal{R}^\mb_{[[\hbar]]})$. The proof (which is true only over
$\mathbb{C}[[\hbar]]$) goes as follows:

Let $x\in H^0(\mathcal{R}^\mb, d_0)$ be a cycle. One needs to
construct a cycle of the form $x+\hbar x_1+\hbar^2 x_2+\dots$ in
$H^0(\mathcal{R}^\mb_{[[\hbar]]}, d_\hbar)$. We a looking for
$$
x^{(k)}=x+\hbar x_1 +\hbar^2 x_2+\dots +\hbar^k x_k
$$
such that

$$
d_\hbar x^{(k)}=0\ mod \ \hbar^{k+1}
$$
Let us note that in this situation one has $d_0((d_\hbar
x^{(k)})_{k+1})=0$. Indeed, it follows from the fact that
$d_\hbar^2(x^{(k)})=0$.

Now we make a step of induction. Consider
$(d_\hbar(x^{(k)}))_{k+1}$. It is $d_0$-cycle, find $x_{k+1}\in
\mathcal{R}^0$ such that $d_0(x_{k+1})=(d_\hbar x^{(k)})_{k+1}$. Set
$$
x^{(k+1)}=x+\hbar x_1+\dots +\hbar^{k+1}x_{k+1}
$$
Then it is clear that $d_\hbar(x^{(k+1)})=0\ mod\ \hbar^{k+2}$.
\end{remark}

\subsection{A proof of the classical Poincar\'{e}-Birkhoff-Witt theorem}
Let $\g$ be a Lie algebra. Its universal enveloping algebra
$\mathcal{U}(\g)$ is defined as the quotient-algebra of the tensor
algebra $T(\g)$ by the two-sided ideal generated by elements
$a\otimes b-b\otimes a-[a,b]$ for any $a,b\in \g$. The
Poincar\'{e}-Birkhoff-Witt theorem says that $\mathcal{U}(\g)$ is
isomorphic to $S(\g)$ as a $\g$-module. We suggest here a (probably
new) proof of this classical theorem, which certainly is not the
simplest one, but sheds some light on the cohomological nature of
the theorem.

Before starting with the proof, let us make some remark. Let us
generalize the universal enveloping algebra as follows. Consider the
tensor algebra $T(x_1,\dots,x_n)$ and its quotient $A_{c_{ij}^k}$ by
the two-sided ideal generated by the relations $x_i\otimes
x_j-x_j\otimes x_i-\sum_{k}c_{ij}^k x_k$, $1\le i<j\le n$, where
$c_{ij}^k$ are not supposed to satisfy the Jacobi identity
\begin{equation}\label{eq2.2.1}
\sum_a(c_{ij}^a c_{ak}^b+c_{jk}^a c_{ai}^b+c_{ki}^a c_{aj}^b)=0
\end{equation}
Then, if (\ref{eq2.2.1}) is not satisfied, the algebra
$A_{c_{ij}^k}$ is smaller than $S(x_1,\dots,x_k)$, that is, the
two-sided ideal, generated by the relations, is bigger than in the
Lie algebra case when (\ref{eq2.2.1}) is satisfied.

Now we pass to the proof. Let $\g$ be a Lie algebra. By the
discussion in Section 1.3, $CoBar^\mb(\Lambda^+(\g))$ is a free
resolution of the symmetric algebra $S(\g)$. Denote the
cobar-differential by $d_0$. Introduce in
$CoBar_\hbar^\mb=CoBar^\mb(\Lambda^+(\g))\otimes\mathbb{C}[\hbar]$ a
new differential $d_0+d_1$, where $d_1\colon CoBar_\hbar^\mb\to\hbar
\cdot CoBar^{\mb+1}_\hbar$ comes from the chain differential in the
Lie homology complex
$\partial\colon\Lambda^i(\g)\to\Lambda^{i-1}(\g)$. We denote
$$
d_1=\hbar\partial
$$
The equation $(d_0+d_1)^2=0$ follows from the fact that the chain
Lie algebra complex is a dg coalgebra, and, therefore, its
cobar-complex is well-defined.

Now, by Lemma 2.1(i), the complex $CoBar^\mb_\hbar$ has only 0
degree cohomology, which is isomorphic to
$H^0(CoBar^\mb(\Lambda^+(\g)))\otimes\mathbb{C}[\hbar]=S(\g)\otimes\mathbb{C}[\hbar]$
as a (filtered) vector space. On the other hand, we can compute 0-th
cohomology of $(CoBar_\hbar^\mb(\Lambda^+(\g)),d_0+d_1)$ directly.
It is the quotient of the tensor algebra
$T(\g)\otimes\mathbb{C}[\hbar]$ by the two-sided ideal generated by
the relations $a\otimes b-b\otimes a-\hbar [a,b]$, $a,b\in\g$.

The specialization of the last isomorphism for $\hbar=1$ gives the
Poincar\'{e}-Birkhoff-Witt theorem.

\begin{remark}
It was a remark of Victor Ginzburg that for a correct
"specialization at $\hbar$=1" we need Lemma 2.1 over polynomials,
that is the case (i) of this Lemma.
\end{remark}

\subsection{}
Consider the following sequence of maps:
\begin{equation}\label{eq2.3.1}
\begin{aligned}
\
&T_\poly(V^*)\xrightarrow{\mathcal{U}_S}\Hoch(S(V))\simeq\Der(CoBar(S(V^*)))\xrightarrow{\Phi_1}\\
&\xrightarrow{\Phi_1}\Der(CoBar(S^+(V^*)))/\Inn(CoBar(S^+(V^*)))
\end{aligned}
\end{equation}
Here the first map is the Kontsevich formality $L_\infty$ morphism
for the algebra $S(V)$, the second isomorphism follows from the
Stasheff's construction, and the third map is the map $\Phi_1$
defined in Section 1.4.

Apply now the composition (\ref{eq2.3.1}) to the vector space $V[1]$
instead of $V^*$.
\begin{lemma}
Let $V$ be a finite-dimensional vector space. Then there is a
canonical isomorphism of the graded Lie algebras $T_\poly(V^*)\simeq
T_\poly(V[1])$.
\begin{proof}
It is straightforward. The map maps $k$-polyvector field with
constant coefficients on $V^*$ to a $k$-linear function on $V[1]$,
and so on.
\end{proof}
\end{lemma}
\begin{remark}
The algebras $S(V)$ and $\Lambda(V^*[-1])$ are Koszul dual, and they
have isomorphic Hochschild comology with all structures (see [Kel]).
\end{remark}
Denote by $\K$ the correspondence $\K\colon T_\poly(V^*)\to
T_\poly(V[1])$ from Lemma. Let $\alpha$ be a polynomial Poisson
bivector on the space $V^*$. By the correspondence from Lemma, we
get a polyvector field $\K(\alpha)$ which in general is not a
bivector, but still satisfies the Maurer-Cartan equation
\begin{equation}\label{eq2.3.3}
[\K(\alpha),\K(\alpha)]=0 \end{equation}

Let us rewrite (\ref{eq2.3.1}) for $V[1]$:
\begin{equation}\label{eq2.3.4}
\begin{aligned}
\
&T_\poly(V[1])\xrightarrow{\mathcal{U}_\Lambda}\Hoch(\Lambda(V^*))\simeq\Der(CoBar(\Lambda(V)))\xrightarrow{\Phi_1}\\
&\xrightarrow{\Phi_1}\Der(CoBar(\Lambda^-(V)))/\Inn(CoBar(\Lambda^-(V)))
\end{aligned}
\end{equation}
 The composition
(\ref{eq2.3.4}) maps the polyvector $\hbar\K(\alpha)$ to a
derivation $d_\hbar$ of degree +1 in
$\Der(CoBar(\Lambda^-(V)))\xxx$, which satisfies the Maurer-Cartan
equation
\begin{equation}\label{eq2.3.2}
(d+d_\hbar)^2=0
\end{equation}
in $\Der/\Inn$, where $d$ is the cobar-differential.

Actually, (\ref{eq2.3.2}) is satisfied in
$\Der(CoBar(\Lambda^-(V)))\xxx$, not only in
$\Der(CoBar(\Lambda^-(V)))\xxx/\Inn(CoBar(\Lambda^-(V)))\xxx$.
Indeed, we suppose that $V$ is placed in degree 0, then
$CoBar(\Lambda^-(V))$ is $\mathbb{Z}_{\le 0}$-graded. Therefore, any
inner derivation has degree $\le 0$, while $d_\hbar$ has degree +1.
We have the following
\subsection{}
\begin{lemma}
Let $\alpha$ be a Poisson bivector on $V^*$, and let $\K(\alpha)$ be
the corresponding Maurer-Cartan polyvector of degree 1 in
$T_\poly(V[1])$. Then (\ref{eq2.3.4}) defines an $\hbar$-linear
derivation $d_\hbar$ of
$CoBar(\Lambda^-(V))\otimes\mathbb{C}[[\hbar]]$ of degree +1
corresponding to $\hbar\K(\alpha)$, such that
$$
(d+d_\hbar)^2=0
$$
where $d$ is the cobar-differential. Moreover, $d_\hbar$ obeys
\begin{equation}\label{eq2.3.5}
d_\hbar(\xi_i\wedge\xi_j)=\hbar\Sym(\alpha_{ij})+\mathcal{O}(\hbar^2)
\end{equation}
where $\xi_i\wedge\xi_j\in \Lambda^{2}(V)$, $\Sym(\alpha_{ij})\in
T(V)$ is the symmetrization, and $\{\xi_i\}$ is the basis in $V[1]$
dual to the basis $\{x_i\}$ in $V^*$ in which
$\alpha=\sum_{ij}\alpha_{ij}\partial_i\wedge\partial_j$
\begin{proof}
We only need to prove (\ref{eq2.3.5}), all other statements are
already proven. We prove it in details in Section 2.7.
\end{proof}
\end{lemma}
\subsection{}
Let $A$ be an $\hbar$-linear associative algebra which is the
quotient of the tensor algebra $T(V)\xxx$ of a vector space $V$ by
the two-sided ideal generated by relations
$$
x\otimes y-y\otimes x=R(x,y)
$$
for any $x,y\in V$, where $R(x,y)\in \hbar T(V)\xxx$. Consider the
following filtration:
\begin{equation}\label{eq2.5.1}
A\supset\hbar A\supset\hbar^2 A\supset\hbar^3 A\supset\dots
\end{equation}
This is clearly an algebra filtration: $(\hbar^k A)\cdot(\hbar^\ell
A)\subset\hbar^{k+\ell}A$. Consider the associated graded algebra
$\mathrm{gr} A$. We say that the algebra $A$ is a
Poincar\'{e}-Birkhoff-Witt (PBW) algebra if $\mathrm{gr}A\simeq
S(V)\xxx$ as a graded $\mathbb{C}[[\hbar]]$-linear algebra.

In general, $\mathrm{gr}A$ is less than
$S(V)\otimes\mathbb{C}[[\hbar]]$, it is a quotient of
$S(V)\otimes\mathbb{C}[[\hbar]]$. One can say that the PBW property
is equivalent to the property that the quotient-algebra has "the
maximal possible size".

\subsection{}
Let $\alpha$ be a polynomial Poisson bivector in $V^*$. In Sections
2.3 and 2.4 we constructed an $\hbar$-linear derivation $d_\hbar$ on
$CoBar(\Lambda^-(V))\xxx$ such that $(d+d_\hbar)^2=0$ where $d$ is
the cobar-differential. By Section 1.3, the cobar-complex
$CoBar(\Lambda^-(V))$ is a free resolution of the algebra $S(V)$, in
particular, the cohomology of $d$ does not vanish only in degree 0
where it is equal to $S(V)$. We are in the situation of Lemma 2.1.
In particular, the dg algebra $(CoBar(\Lambda^-(V))\xxx,d+d_\hbar)$
has only non-vanishing cohomology in degree 0, and this 0-degree
cohomology is an algebra, which is a PBW algebra by Lemma 2.1.
\begin{theorem*}
The construction above constructs from a Poisson polynomial bivector
$\alpha$ on $V^*$ an algebra $A_\alpha$ with generators
$x_1,\dots,x_n$ and relations $[x_i,x_j]=d_\hbar(\xi_i\wedge\xi_j)$.
This algebra is a PBW algebra. \qed
\end{theorem*}

\begin{conjecture}
The algebra $A_\alpha$ is isomorphic to the Kontsevich star-algebra
on $S(V)\xxx$ constructed from the Poisson bivector $\alpha$. (We
suppose that in the formality morphisms $\mathcal{U}_\Lambda\colon
T_\poly(V[1])\to\Hoch(\Lambda(V^*))$ in (\ref{eq2.3.4}), and
$\mathcal{U}_S\colon T_\poly(V^*)\to\Hoch(S(V))$ which is used in
the construction of the star-product, one uses the same propagator
in the definition of the Kontsevich integrals, see [K97]). The
isomorphism is given by the symmetrization map.
\end{conjecture}

\subsection{An explicit formula}
One can write down explicitly the relations in the algebra
$A_\alpha$, in the terms of the Kontsevich integrals [K97]. For this
we need to find explicit formula for the $\hbar$-linear derivation
$d_\hbar$ in $CoBar(\Lambda^-(V))\xxx$. Here we suppose some
familiarity with [K97].

First of all, recall how the Kontsevich deformation quantization
formula is written. Let $\alpha$ be a Poisson structure on $V^*$.
Then the formula is
\begin{equation}\label{eq2.7.1}
f\star g=f\cdot g+\sum_{k\ge 1}\hbar^k(\sum_{m\ge 1}\frac
1{m!}\sum_{\Gamma \in G_{2,m}^2}W_\Gamma
U_\Gamma(\alpha,\dots,\alpha))
\end{equation}
Here $\Gamma$ is an admissible graph with two vertices on the "real
line" and $m$ vertices in the upper half-plane, and with 2 outtgoing
edges at each vertex in the upper half-plane, that is, $\Gamma\in
G_{2,m}^2$, in particular, it is an oriented graph with $2m$ edges;
$W_\Gamma$ is the Kontsevich integral of the graph $\Gamma$. {\it
Let us note that the all graphs involved in (\ref{eq2.7.1}) may have
arbitrary many incoming edges at each vertex at the upper
half-plane, and exactly two outgoing edges}.

Now let $\alpha=\sum_{ij}\alpha_{ij}\partial_i\wedge \partial_j$,
where $\alpha_{ij}=\sum_I c_{ij}^I x_{i_1}\dots x_{i_k}$ ($I$ is a
multi-index).

Then the "Koszul dual" polyvector $\K(\alpha)$ is a polyvector field
with quadratic coefficients:

\begin{equation}\label{eq2.7.2}
\K(\alpha)=\sum_{i,j,I}c_{ij}^I (\xi_i\xi_j) \cdot
\partial_{\xi_{i_1}}\wedge\dots\wedge\partial_{\xi_{i_k}}
\end{equation}
It has total degree 1 and satisfies the Maurer-Cartan equation.

Firstly we write the formula for the image of $\K(\hbar\alpha)$ by
the Kontsevich formality, that is, denote by
\begin{equation}\label{eq2.7.3}
\mathcal{U}(\K(\alpha))=\hbar U_1(\K(\alpha))+\hbar^2
\frac12U_2(\K(\alpha),\K(\alpha))+\dots+\hbar^k\frac1{k!}U_k(\K(\alpha),\dots,\K(\alpha))+\dots
\end{equation}
We can write down explicitly this formula in graphs. {\it Let us
note the the graphs involving in (\ref{eq2.7.3}) may have arbitrary
many outgoing edges in the vertices at the upper half-plane, but
exactly two incoming edges, because all components of $\K(\alpha)$
are quadratic polyvector fiels}, and by a simple dimension
computation. That is, in a sense the graphs in (\ref{eq2.7.3}) are
"dual" to the graphs in (\ref{eq2.7.1}).

Let us note also that the right-hand side of (\ref{eq2.7.3}) is a
polydifferential operator in $\Hoch(\Lambda(V^*))$ of
non-homogeneous Hochschild degree, but of the total (Hochschild
degree and $\Lambda$-degree) +1.

Now we should apply to $\mathcal{U}(\K)$ our map $\Phi_1$ to get a
derivation of the cobar-complex $CoBar(\Lambda^-(V))$. After this,
we get the final answer for $d_\hbar$.

Let us compute its component in the first power of $\hbar$. It is
just the Hochschild-Kostant-Rosenberg map of $U_1(\K(\alpha))$ which
is the symmetrization map $\Sym\colon S(V)\to T(V)$ in this case.

Let us note that in the case of a quadratic Poisson structure
algebra the relations $R_{ij}\in (V\otimes V)\xxx\subset T(V)\xxx$
are quadratic.

Faculty of Science, Technology and Communication, Campus
Limpertsberg, University of Luxembourg,
162A avenue de la Faiencerie, L-1511 LUXEMBOURG\\
{\it e-mail}: {\tt borya$\_$port@yahoo.com}

\end{document}